\documentclass[12pt]{article}      

\usepackage[centertags]{amsmath}
\usepackage{amsfonts}
\usepackage{amssymb}

\renewcommand{\sectionmark}[1]
                    {\markboth{Kapitel \thesection\ #1}{}}
\renewcommand{\sectionmark}[1]
                 {\markright{} }

\newtheorem{thm}{Theorem}
\newtheorem{lem}[thm]{Lemma}
\newtheorem{corollary}[thm]{Corollary}

\newtheorem{fremdersatz}{Theorem}

\newenvironment{satOrig}[1]
        {\pagebreak[2] \begin{fremdersatz} {\bf #1} \quad\sl}
        { \end{fremdersatz}}

\numberwithin{equation}{section}

\newenvironment{proofof}[1]
        {\pagebreak[2] \vspace{-1pt}{\bf Proof#1.}  }
        {\hfill $\blacksquare$ \vspace{2pt}}

\newcommand{\fa}{\mathcal{F}}
\newcommand{\la}{\mathcal{L}}
\newcommand{\pa}{\mathcal{P}}
\newcommand{\nat}{{\rm I\! N}}

\newcommand{\co}{{\mathbb C}}

\newcommand{\gl}{\left\{}
\newcommand{\gr}{\right\}}
\newcommand{\kl}{\left(}
\newcommand{\kr}{\right)}
\newcommand{\kj}{\overline}
\renewcommand{\l}{\left}
\renewcommand{\r}{\right}
\newcommand{\limn}{\lim_{n\to\infty}}

\newcommand{\In}{\subseteq}
\newcommand{\mi}{\setminus}
\newcommand{\abb}{\longrightarrow}

\newcommand{\beq}{\begin{equation}}
\newcommand{\eeq}{\end{equation}}
\newcommand{\beqar}{\begin{eqnarray}}
\newcommand{\eeqar}{\end{eqnarray}}
\newcommand{\beqaro}{\begin{eqnarray*}}
\newcommand{\eeqaro}{\end{eqnarray*}}
\newcommand{\bsat}{\begin{thm}}
\newcommand{\esat}{\end{thm}}
\newcommand{\bsatorig}{\begin{satOrig}}
\newcommand{\esatorig}{\end{satOrig}}
\newcommand{\blem}{\begin{lem}}
\newcommand{\elem}{\end{lem}}
\newcommand{\bkor}{\begin{corollary}}
\newcommand{\ekor}{\end{corollary}}
\newcommand{\bbew}{\begin{proofof}}
\newcommand{\ebew}{\end{proofof}}

\renewcommand{\rho}{\varrho}
\renewcommand{\phi}{\varphi}
\renewcommand{\epsilon}{\varepsilon}

\begin{document}

\title{Differential inequalities and quasi-normal families}

\author{Roi Bar, J\"urgen Grahl and Shahar Nevo}

\date{\today}

\maketitle

  \begin{center}
{\it     In honor of Professor Lawrence Zalcman on the occasion
  of his 70th birthday  }
\end{center}

\begin{abstract}
We show that a family $\fa$ of meromorphic functions in a domain $D$
satisfying 
$$\frac{|f^{(k)}|}{1+|f^{(j)}|^\alpha}(z)\ge C \qquad \mbox{ for all }
z\in D \mbox{ and all } f\in\fa$$ 
(where $k$ and $j$ are integers with $k>j\ge 0$ and $C>0$, $\alpha>1$
are real numbers) is quasi-normal. Furthermore, if all functions in
$\fa$ are holomorphic, the order of quasi-normality of $\fa$ is at most $j-1$. 
The proof relies on the Zalcman rescaling method and previous results
on differential inequalities constituting normality. 


{\bf Mathematics Suject Classification:} 30D45, 30A10
\end{abstract}

\section{Introduction and statement of results}

According to Marty's theorem, a family $\fa$ of meromorphic functions in a
domain $D\In\co$ is normal (in the sense of Montel) if and only if the
family $\fa^\#:=\gl f^\#\;:\;f\in\fa\gr$ of the corresponding
spherical derivatives is locally uniformly bounded in $D$; here,
$f^\#$ is defined by $f^\#:=\frac{|f'|}{1+|f|^2}$.

In \cite{GrahlNevo-Spherical} we studied families of meromorphic
functions whose spherical derivatives are bounded away from zero and
proved the following counterpart to Marty's theorem. 

\bsatorig{}\label{MartyBelow}
Let $D\In\co$ be a domain and $C>0$. Let $\fa$ be a family
of functions meromorphic in $D$ such that 
$$f^\#(z)\ge C \qquad\mbox{ for all } z\in D \mbox{ and all }
f\in \fa.$$
Then $\fa$ is normal in $D$. 
\esatorig

Hence, the condition $\frac{|f'|}{1+|f|^2}(z)=f^\#(z)\ge C$ can be considered
as a differential inequality that constitutes normality. In
\cite{LiuNevoPang}, \cite{ChenNevoPang}, \cite{GNP-NonExplicit} and
\cite{GN-Marty} we studied more general
differential inequalities, involving higher derivatives, with respect
to the question whether they constitute normality or at least
quasi-normality. 

Before summarizing the main results from these studies, as far as they
are relevant in the context of the present paper, we would like to remind
the reader of the definition of quasi-normality and also to introduce
some notations.  A family $\fa$ of meromorphic functions in a domain
$D\In\co$ is said to be {\it quasi-normal} if from each sequence $\gl
f_n\gr_n$ in $\fa$ one can extract a subsequence which converges
locally uniformly (with respect to the spherical metric) on $D\mi E$
where the set $E$ (which may depend on $\gl f_n\gr_n$) has no
accumulation point in $D$. If the exceptional set $E$ can always be
chosen to have at most $q$ elements, we say that $\fa$ is {\it
  quasi-normal of order at most} $q$. Finally, $\fa$ is said to be
{\it quasi-normal of (exact) order $q$} if it is quasi-normal of order
at most $q$, but not quasi-normal of order at most $q-1$.

C.T. Chuang has extended the concept of quasi-normality by introducing
the notion of $Q_m$-normal families (\cite{Chuang}, see also
\cite{Nevo-Qm}). A family $\fa$ of meromorphic functions in a
domain $D$ is called $Q_m$-normal if from each sequence in $\fa$ one
can extract a subsequence which converges locally uniformly (with
respect to the spherical metric) on $D\mi E$ where the set $E$
satisfies $E^{(m)}\cap D=\emptyset$; here $E^{(1)}=E'$ is the derived
set of $E$, i.e. the set of its accumulation points, and $E^{(m)}$ is
defined inductively by $E^{(m)}:=(E^{(m-1)})'$ for $m\ge 2$. A
$Q_1$-normal family is just a quasi-normal family.

For $z_0\in\mathbb C$ and $r>0,$ we set
$\Delta(z_0,r):=\{z\in\co:|z-z_0|<r\} $ and
$\Delta'(z_0,r):=\Delta(z_0,r)\mi\gl z_0\gr. $ Furthermore, we denote
the open unit disk by $\Delta:=\Delta(0,1)$. 
We write ``$f_n\Longrightarrow f$ on $D$''
to indicate that the sequence $\{f_n\}_n$ converges to $f$ uniformly
on compact subsets of $D$  (w.r.t. the Euclidean metric).

Now we can turn to the results on differential inequalities and
quasi-normality known so far. While \cite{GNP-NonExplicit} and
\cite{GN-Marty} dealt with generalizations of Marty's theorem (more
precisely with conditions of the form
$\frac{|f^{(k)}|}{1+|f|^\alpha}(z)\le C$), in \cite{ChenNevoPang} the
following extension of Theorem \ref{MartyBelow} was proved.

\bsatorig{\cite{ChenNevoPang}} \label{CNP}
Let $\alpha> 1$ and $C>0$ be real numbers and $k\ge 1$ be an
integer. Let $\fa$ be a family of meromorphic functions in some domain
$D$ such that  
\beq\label{LowerEstimate}
\frac{|f^{(k)}|}{1+|f|^\alpha}(z)\ge C \qquad \mbox{ for all }
z\in D \mbox{ and all } f\in\fa.
\eeq
Then $\fa$ is normal.
\esatorig

This result doesn't hold any longer if $\alpha>1$ is replaced by
$\alpha=1$ as the family $\gl z\mapsto nz^k\;: \; n\in\nat\gr$ which
is not normal at $z=0$ demonstrates. However, at least for $k=1$ 
condition (\ref{LowerEstimate}) implies quasi-normality if
$\alpha=1$, as shown in \cite{LiuNevoPang}. 

\bsatorig{\cite{LiuNevoPang}} \label{LNP}
Let $C>0$ be a real number and $\fa$ be a family of meromorphic functions in some domain
$D$ such that  
$$
\frac{|f'|}{1+|f|}(z)\ge C \qquad \mbox{ for all }
z\in D \mbox{ and all } f\in\fa.
$$
Then $\fa$ is quasi-normal.
\esatorig
 
In the present paper we prove the following extension of Theorem \ref{CNP}. 

\bsat{}\label{mainresult}
Let $k$ and $j$ be integers with $k>j\ge 0$ and $C>0$, $\alpha>1$
be real numbers. Let $\fa$ be a family of
meromorphic functions in some domain $D$ such that for each $f\in\fa$
\beq\label{LowerBound}
\frac{|f^{(k)}|}{1+|f^{(j)}|^\alpha}(z)\ge C \qquad \mbox{ for all }
z\in D.
\eeq
Then $\fa$ is quasi-normal in $D$. 

If all functions in $\fa$ are holomorphic, $\fa$ is quasi-normal of
order at most $j-1$. (For $j=0$ and $j=1$ this means that it is normal.)
\esat

To simplify notations, for an arbitrary family of holomorphic or
meromorphic functions in a domain $D$, integers $k>j\ge
0$ and a real number $\alpha\ge0$ we define
$$\fa_{k,j,\alpha}:=\gl \frac{|f^{(k)}|}{1+|f^{(j)}|^\alpha} : f\in\fa\gr.$$
Then Theorem \ref{mainresult} states that a family $\fa$ of
meromorphic functions is quasi-normal if $\fa_{k,j,\alpha}$ (where
$k>j\ge 0$, $\alpha>1$) is uniformly bounded away from zero. (Of
course, since quasi-normality is a local property, it also suffices to
assume that $\fa_{k,j,\alpha}$ is {\it locally} uniformly bounded away
from zero.)

We note that the ``meromorphic case'' (more precisely: the
non-holomorphic case) of Theorem \ref{mainresult} deviates from the
holomorphic case only if $k\ge\alpha (j+1)-1$. Indeed, if $f$ is a
meromorphic function in $D$ satisfying (\ref{LowerBound}) and if $f$
has a pole of order $p\ge 1$ 
at $z_0\in D$, then $f^{(j)}$ has a pole of order $p+j$ and $f^{(k)}$
has a pole of order $p+k$ at $z_0$. From (\ref{LowerBound}) we see that
\beq\label{MeroNontrivial} p+k-\alpha(p+j)\ge 0, \quad \mbox{ i.e. }
\quad 1\le p \le \frac{k-\alpha j}{\alpha-1}, \eeq which implies
$k\ge\alpha (j+1)-1$. In other words, if $k<\alpha (j+1)-1$, the
functions satisfying (\ref{LowerBound}) cannot have any poles. 

Concerning the order of quasi-normality, Theorem \ref{mainresult} is
sharp, both in the holomorphic and in the meromorphic case: 
\begin{itemize}
\item[(1)] 
For any given $a\in\co\mi\gl0\gr$ there exists a $C>0$ such that all functions
$f(z):=p(z)+e^{az}$, where $p$ is an arbitrary polynomial of degree
$j-1$, satisfy (\ref{LowerBound}) in the unit disk $\Delta$. Since the class of
polynomials of degree $j-1$ is quasi-normal of order exactly $j-1$
(Lemma \ref{PolynQuasi}), the order of quasi-normality of the class of
holomorphic functions in $\Delta$ satisfying (\ref{LowerBound}) cannot
be less than $j-1$.  

So in the holomorphic case Theorem \ref{mainresult} provides a
differential inequality that separates between different orders of
quasi-normality. 
\item[(2)] 
As to the meromorphic case, in the situation of Theorem
\ref{mainresult} the order of quasi-normality can be arbitrarily
large: Indeed, for given $m\in\nat$, consider the  functions 
$$f_n(z)=\frac{1}{z-z_1}+\dots+\frac{1}{z-z_m}+n,$$
where
$z_1,\dots,z_m$ are distinct points in $\Delta$. Then
$f_n^{(k)}=f_1^{(k)}$ has only finitely many zeros $w_1,\dots,w_\ell$,
and they are independent of $n$. Thus if we take $D:=\Delta\mi
\bigcup_{\mu=1}^{\ell} \kj{\Delta(w_\mu,\delta)}$ for small enough
$\delta>0$, then all functions $f_n$ satisfy (\ref{LowerBound}) for
an appropriate $C>0$, provided that $k\ge \alpha(j+1)-1$. Obviously,
the (only) points of non-normality of every subsequence of $\gl
f_n\gr_n$ are the points $z_1,\dots,z_m$. So $\gl f_n\gr_n$ is a
quasi-normal sequence of order $m$. 
\end{itemize}

Furthermore, at least for $j \ge 1$ the condition $\alpha>1$ in
Theorem \ref{mainresult} cannot be replaced by $\alpha\ge 1$. This is
shown by the sequence of the functions $f_n(z):=z^n-3^n$ on the
annulus $D:=\gl z\in\co\;:\; 2<|z|<4\gr$. Obviously, 
$\gl\frac{|f_n^{(k)}|}{1+|f_n^{(j)}|}\gr_n$ tends to $\infty$ locally
uniformly on $D$ whenever $k>j\ge 1$, but $\gl f_n\gr_n$ is not normal
at any point $z$ with $|z|=3$, hence it isn't quasi-normal (and not
even $Q_m$-normal for any $m$). 

We conjecture that the first part of Theorem \ref{mainresult}
(concerning the quasi-normality of $\fa$, but without the statement
about the order of quasi-normality) remains valid for $j=0$ and
$\alpha=1$, but we were not able to prove this in general. (For $k=1$
this conjecture is true by Theorem \ref{LNP}.)

Surprisingly, it turns out that the condition
$\frac{|f^{(k)}|}{1+|f|^\alpha}(z)\ge C$ for all $z\in D$ (where
$k\ge 1$, $\alpha\ge 1$) cannot hold for all functions $f$ in a
certain (infinite) family if there is some point $z_0\in D$ such that $\fa$ is
``not holomorphic'' and ``not zero-free'' at $z_0$, in the sense that
in arbitrary small neighborhoods of $z_0$ all but finitely many
functions in $\fa$ have both zeros and poles. More precisely we have
the following result.

\bsat{} \label{Phenomenon}
Let $\alpha\ge1$ be a real number. Let $\gl f_n\gr_n$ be a sequence of
meromorphic functions in a domain $D$ and let some $z_0\in D$ be
given. Assume that there exist sequences $\gl z_n\gr_n$ and $\gl
p_n\gr_n$ in $D$ such that $\limn z_n=\limn p_n=z_0$ and $f_n(z_n)=0$,
$f_n(p_n)=\infty$ for all $n$. Then  
$$\inf_{z\in D} \frac{|f_n'|}{1+|f_n|^\alpha}(z)\abb0 \quad (n\to\infty).$$
\esat

From this observation and Gu's normality criterion (Lemma \ref{Gu}) we
easily obtain the following result. 

\bkor\label{Phenomenon-Cor}
Let $\alpha\ge1$ and $C>0$ be real numbers. Let $\gl f_n\gr_n$ be a sequence of
meromorphic functions in a domain $D$ satisfying 
$$\frac{|f_n'|}{1+|f_n|^\alpha}(z)\ge C\qquad\mbox{ for all } z\in D
\mbox{ and all } n\in\nat.$$ 
If there exists a sequence $\gl p_n\gr_n$ in $D$ such that $\limn
p_n=z_0$ and  $f_n(p_n)=\infty$ for all $n$, then  $\fa$ is normal at $z_0$. 
\ekor

\section{Some Lemmas}

The most essential tool in our proofs is a famous rescaling lemma
which was originally proved by L.~Zalc\-man~\cite{zalcman} and later
extended by X.-C.~Pang~\cite{pang89,pang90} and by H.~Chen and Y.~Gu
\cite{chengu}. Here we require the following version from
\cite{pangzalc2000a} (see also \cite{GrahlNevo-Spherical} for a proof
of the direction ``$\Leftarrow$'').

\blem[Zalcman-Pang Lemma] \label{zalclemma} 
Let $\fa$ be a family of meromorphic functions in a domain $D$ all of
whose zeros have multiplicity at least $m$ and all of whose poles have
multiplicity at least $p$ and let $-p<\beta<m$. Then $\fa$ is not
normal at some $z_0\in D$ if and only if there exist sequences $\gl
f_n\gr_n\In \fa$, $\gl z_n\gr_n \In D$ and $\gl \rho_n\gr_n \In (0,1)$
such that $\gl\rho_n\gr_n$ tends to 0, $\gl z_n\gr_n$ tends to $z_0$ and such
that the sequence $\gl g_n\gr_n$ defined by
$$g_n(\zeta):= \frac{1}{\rho_n^\beta}\cdot f_n(z_n+\rho_n \zeta)$$
converges locally uniformly in $\co$ (with respect to the spherical
metric) to a non-constant function $g$ meromorphic in $\co$.
\elem

Furthermore, we need the following famous normality criterion due to Y. Gu \cite{Gu}.

\blem\label{Gu} 
Let $k\ge 1$ be an integer. Then the family of all functions $f$
meromorphic in a domain $D\In \co$ satisfying $f(z)\ne 0$,
$f^{(k)}(z)\ne 1$ for every $z\in D$ is normal. 
\elem

\blem{}\label{PolynQuasi}
For every integer $k\ge 0$ the family $\pa_k$ of polynomials of degree at most
$k$ is quasi-normal of exact order $k$ in every domain in $\co$. 
\elem

\bbew{}
This is a special case of Theorem A.5 in \cite{Schiff} which states
that any family of holomorphic functions in a domain $D$ which do not
take a value $a\in\co$ more than $p$ times nor a value $b\in\co$ more
than $q$ times is quasi-normal of order at most $\min\gl p,q\gr$. That
the order of quasi-normality of $\pa_k$ cannot be less than $k$
follows by considering the sequence $\gl n\cdot p\gr_n$ where $p(z):=(z-z_1)\cdot\ldots\cdot
(z-z_k)$ with arbitrarily chosen distinct points $z_1,\dots,z_k\in
D$.  
\ebew

\section{Proofs}

\bbew{ of Theorem \ref{mainresult}}
{\bf I.} We start with the proof of the holomorphic case, i.e. we show
that a family $\fa$ of holomorphic functions in some domain $D$
satisfying (\ref{LowerBound}) is quasi-normal of order at most $j-1$. 

{\bf (a)} For $j=0$ this is true by Theorem \ref{CNP}. (In fact, in
this case we have normality for any $k\ge 1$.)

Now we consider the case $j=1$ (which turns out to require the main
efforts of the proof of the holomorphic case). Here we have $k\ge 2$. 

We assume that $\fa$ is not normal at some
$z_0\in D$. Then we choose some $\beta>0$ such that $(\beta+1)\cdot
\alpha-(\beta+k)>0$ and find by the Zalcman-Pang lemma (Lemma \ref{zalclemma})
sequences $\gl f_n\gr_n\In \fa$, $\gl z_n\gr_n \In D$ and
$\gl \rho_n\gr_n \In (0,1)$ such that $\gl\rho_n\gr_n$ tends to 0,
$\gl z_n\gr_n$ tends to $z_0$ and such that the sequence $\gl
g_n\gr_n$ defined by
$$g_n(\zeta):= \rho_n^\beta\cdot f_n(z_n+\rho_n \zeta)$$
converges locally uniformly in $\co$ (with respect to the spherical
metric) to a non-constant function $g$ meromorphic in $\co$. (This
choice of $\beta$ is admissible since the $f_n$ have no poles.)
Differentiating gives 
$$\rho_n^{\beta+1}\cdot f_n'(z_n+\rho_n \zeta)\Longrightarrow
g'(\zeta) \qquad (n\to\infty).$$
Now by Theorem \ref{CNP} the derivatives $f_n'$ form a normal family
(since $f_n'$ satisfies (\ref{LowerEstimate}) with $k-1$ instead of
$k$). Therefore by the converse of the Zalcman-Pang lemma we deduce
that $g'$ is constant. So $g$ has the 
form $g(\zeta)=A\zeta+B$ where $A\ne 0$. Differentiating $k-1(\ge 1)$
more times gives 
$$\rho_n^{\beta+k}\cdot f_n^{(k)}(z_n+\rho_n \zeta)\Longrightarrow0 
\qquad (n\to\infty).$$
We fix some $\zeta_0\in\co$. Then for $n$ large enough we have 
$$\rho_n^{\beta+k}\cdot |f_n^{(k)}(z_n+\rho_n \zeta_0)|\le 1 
\qquad \mbox{ and } \qquad 
\rho_n^{\beta+1}\cdot|f_n'(z_n+\rho_n \zeta_0)|\ge \frac{|A|}{2},$$
hence
$$\frac{|f_n^{(k)}|}{1+|f_n'|^\alpha}(z_n+\rho_n\zeta_0)
\le\kl\frac{2}{|A|}\kr^\alpha\cdot \rho_n^{(\beta+1)\alpha-(\beta+k)}.$$
Here, the condition $(\beta+1)\cdot \alpha-(\beta+k)>0$ ensures that
the right hand side tends to 0 for $n\to\infty$, contradicting (\ref{LowerBound}).  

Hence $\fa$ is normal if $j=1$ and all functions in $\fa$ are holomorphic. 

{\bf (b)} We prove the following claim.

{\bf Claim:} If $j\ge 1$, every function $f\in \fa$ can be written as
$f=h_f+p_f$ where $p_f$ is 
a polynomial of degree at most $j-1$ and $\gl h_f \;:\;
f\in\fa\gr$ is a locally uniformly bounded family of holomorphic
functions. 

{\bf Proof.} We fix an arbitrary $z_0\in D$. For given $f\in \fa$, we
can find a polynomial $p_f$ of degree at most $j-1$ such that
$p_f^{(\mu)}(z_0)=f^{(\mu)}(z_0)$ for $\mu=0,\dots,j-1$ (i.e. $p_f$ is
the $(j-1)$'th Taylor polynomial of $f$ at $z_0$). Then for
$h_f:=f-p_f$ we have $h_f^{(\mu)}(z_0)=0$ for $\mu=0,\dots,j-1$.
Furthermore, $p_f^{(j)}\equiv 0$ implies $h_f^{(j)}=f^{(j)}$ and
$h_f^{(k)}=f^{(k)}$, so $h_f$ satisfies (\ref{LowerBound}) as well.
But this means that $h_f^{(j-1)}=:\varphi$ satisfies
$\frac{|\varphi^{(k-j+1)}|}{1+|\varphi'|^\alpha}(z)\ge C$ for all $z\in D$, a
condition which implies normality as we have shown in (a).  Hence $\gl
h_f^{(j-1)}\;:\; f\in\fa\gr$ is normal. Observing the normalization
$h_f^{(j-1)}(z_0)=0$ we obtain that $\gl h_f^{(j-1)}\;:\; f\in \fa\gr$
is even locally uniformly bounded. Using $h_f^{(\mu)}(z_0)=0$ for
$\mu=0,\dots,j-2$, we inductively deduce that $\gl h_f^{(j-2)}\;:\;
f\in \fa\gr, \dots, \gl h_f\;:\; f\in \fa\gr$ are uniformly bounded.

{\bf (c)} Now we can complete the proof of the holomorphic case for
$j\ge 2$. Let some sequence $\la$ in $\fa$ be given. Then by the Claim
in (b) and by Montel's theorem there is a subsequence $\gl f_n\gr_n$
such that each $f_n$ can be written in the form $f_n=h_n+p_n$ where
$p_n$ is a polynomial of degree at most $j-1$, $h_n$ is holomorphic in
$D$ and the sequence $(h_n)_n$ converges locally uniformly in $D$ to
some holomorphic limit function $h$. Now Lemma \ref{PolynQuasi}
immediately implies that $\gl f_n\gr_n$ is quasi-normal of order at
most $j-1$.  This shows the assertion.

{\bf II.} Now we turn to the general case. Let some point $z_0\in D$
and some sequence $\la$ in $\fa$ be given. If
there is a $\delta>0$ such that all but finitely many functions in
$\la$ are holomorphic in $\Delta(z_0,\delta)$, then $\la$ is quasi-normal
at $z_0$ by the holomorphic case~{\bf I.} 

Otherwise we can find a subsequence $\gl f_n\gr_n$ of $\la$ and a
sequence $\gl z_n\gr_n$ converging to $z_0$ such that $z_n$ is a pole
of $f_n$ for each $n$. Let $p_n\ge 1$ be the multiplicity of $z_n$ as a
pole of $f_n$. Then as in (\ref{MeroNontrivial}) we see 
$$p_n \le \frac{k-\alpha j}{\alpha-1}$$
for all $n$. In particular, the sequence $\gl p_n\gr_n$ is
bounded. After turning to an appropriate subsequence we may assume
that it is constant, $p_n=p$ for all $n$.  

Using the estimate $x\le 1+x^\alpha$ for all $x\ge 0$ we obtain 
$$\l|\frac{f_n^{(j)}}{f_n^{(k)}}(z)\r|
\le \frac{1+|f_n^{(j)}|^\alpha}{|f_n^{(k)}|}(z)\le \frac{1}{C}
\qquad\mbox{ for all } z\in D, n\in\nat.$$
So by Montel's theorem the sequence
$\gl\frac{f_n^{(j)}}{f_n^{(k)}}\gr_n$ is normal, and w.l.o.g. we may
assume that it tends to some holomorphic limit function $h$ locally
uniformly in $D$. Obviously, $z_n$ is a zero of
$\frac{f_n^{(j)}}{f_n^{(k)}}$ of order $k-j$, and near $z_n$ we have
the Taylor expansion
$$\frac{f_n^{(j)}}{f_n^{(k)}}(z)
=\frac{(-1)^{k-j}}{(p+j)\cdot(p+j+1)\cdot\ldots\cdot(p+k-1)}\cdot
(z-z_n)^{k-j}\cdot (1+O(z-z_n)),$$
which yields 
$$\kl\frac{f_n^{(j)}}{f_n^{(k)}}\kr^{(k-j)}(z_n)=\frac{(-1)^{k-j}\cdot(k-j)!}{(p+j)\cdot(p+j+1)\cdot\ldots\cdot(p+k-1)}.$$
Since this quantity does not depend on $n$, we deduce
$h^{(k-j)}(z_0)\ne0$. So $h$ is not constant and has a zero of
multiplicity $k-j$ at $z_0$. This ensures that $\gl p_n\gr_n$ is the
only sequence of poles of the functions $f_n$ that tends to $z_0$. 

We fix some $r>0$ such that $h$ has no zeros in $\Delta(z_0,r)\mi\gl
z_0\gr$. Then for any $\delta_1,\delta_2$ with $0<\delta_1<\delta_2<r$
we conclude that for all $n$ sufficiently large $f_n$ is holomorphic
in the annulus $\gl z\in\co\;:\; \delta_1<|z-z_0|<\delta_2\gr$. By
applying the holomorphic case {\bf I.} we deduce that $\gl f_n\gr_n$
is quasi-normal of order at most $j-1$ in each such annulus. But this
means that it is quasi-normal of order at most $j-1$ in the punctured
disk $\Delta'(z_0,r)$, and so it is quasi-normal of order at most $j$
in $\Delta(z_0,r)$. Here it is crucial to have information about the
order of quasi-normality in $\Delta'(z_0,r)$.  Otherwise we could only
deduce that $\gl f_n\gr_n$ is $Q_2$-normal in $\Delta(z_0,r)$.

Since this holds for any $z_0\in D$ and any sequence $\la$ in
$\fa$, we have shown that  $\fa$ is quasi-normal. 
\ebew

\bbew{ of Theorem \ref{Phenomenon}}
We assume that the assertion is wrong. Then after turning to an
appropriate subsequence we may assume that there is some
$\varepsilon>0$ such that
$$\frac{|f_n'|}{1+|f_n|^\alpha}(z) \ge \varepsilon \qquad\mbox{ for
  all } z\in D \mbox{ and all } n\in\nat,$$
hence
$$\l|\frac{f_n'}{f_n}(z)\r| \ge \varepsilon \qquad\mbox{ for
  all } z\in D \mbox{ and all } n\in\nat.$$
But this means that $\gl \frac{f_n}{f_n'}\gr_n$ is uniformly
  bounded, hence normal. W.l.o.g. we may assume
  $\frac{f_n}{f_n'}\Longrightarrow H$ ($n\to\infty$) for some
  holomorphic function $H$. For each $n$, both $p_n$ and $z_n$ are
  zeros of $\frac{f_n}{f_n'}$. Thus, if $H\not\equiv 0$, from
  Hurwitz's theorem we obtain that $z_0$ is a zero of $H$ of
  multiplicity at least 2. Hence we have $H'(z_0)=0$, and this also
  holds if $H\equiv 0$. On the other hand, in view of $f_n(z_n)=0$ and
  $|f_n'(z_n)|\ge \varepsilon$ we have 
$$\kl\frac{f_n}{f_n'}\kr'(z_n)=1-\frac{f_nf_n''}{(f_n')^2}(z_n)=1$$
for all $n$, hence $H'(z_0)=\limn \kl\frac{f_n}{f_n'}\kr'(z_n)=1$ by
Weierstra\ss{}'s theorem. This is a contradiction.
\ebew

\bbew{ of Corollary \ref{Phenomenon-Cor}}
By Theorem \ref{Phenomenon} there exists a neighborhood
$\Delta(z_0,\delta)$ of $z_0$ such that all but finitely many $f_n$
are zero-free in $\Delta(z_0,\delta)$. This and the fact that
$f_n'(z)\ne \frac{C}{2}$ for all $z\in D$ and all $n$ yields the
normality of $\gl f_n\gr_n$ in $\Delta(z_0,\delta)$ by Gu's normality criterion
(Lemma \ref{Gu}).  
\ebew

\end{document}